\documentclass[12pt]{article}
\usepackage{fullpage}
\usepackage{epstopdf}
\usepackage{amsmath,amssymb,graphicx,float,epsf}
\usepackage{epic,eepic}
\usepackage{float}
\restylefloat{figure}

\newcommand{\rf}[1]{(\ref{#1})}
\newcommand{\beq}[1]{ \begin{equation}\label{#1} }
\newcommand{\eeq}{\end{equation}}
\def\NN{\hbox{I\kern-.2em\hbox{N}}}
\def\RR{\hbox{I\kern-.2em\hbox{R}}}
\newcommand{\qed}{\hbox to 0pt{}\hfill$\rlap{$\sqcap$}\sqcup$\vspace{2mm}}
\newtheorem{uess}{Lemma}
\newtheorem{guess}{Theorem}
\newtheorem{remark}{Remark}
\newtheorem{corol}{Corollary}

\title{Absolute and Delay-Dependent Stability of Equations with a 
Distributed Delay: a Bridge from Nonlinear Differential to Difference 
Equations}
\author{Elena Braverman$^{1,2}$ and Sergey Zhukovskiy$^1$ \\ 
Department of Mathematics and Statistics, University of Calgary, \\
2500 University Drive N.W., Calgary, AB, Canada T2N 1N4
}
\date{}

\begin{document}
\maketitle

\footnotetext[1]{Partially supported by the NSERC Research Grant}
\footnotetext[2]{Corresponding author. E-mail {\em
maelena@math.ucalgary.ca}. Fax (403)-282-5150. Phone (403)-220-3956.}

\begin{abstract}
We study delay-independent stability in nonlinear models with a 
distributed delay which have a positive equilibrium.
Such models frequently occur in population dynamics and other 
applications. In particular, we construct a relevant difference 
equation such that its stability implies stability
of the equation with a distributed delay and a finite memory.
This result is, generally speaking, incorrect for systems with infinite
memory. If the relevant difference equation is unstable, we describe 
the general delay-independent attracting set and also demonstrate that
the equation with a distributed delay is stable
for small enough delays.  
%Instability of a positive equilibrium of the difference equation
%does not imply that 
\end{abstract}
  


\section{Introduction}

In models of population dynamics which are described by an autonomous
differential equation 
\begin{equation}
\label{01}
\frac{dN}{dt}=f(N)-g(N),
\end{equation}
where $f(N)$ and $g(N)$ are reproduction and mortality rates, 
respectively, $f(N)>0$, $g(N)>0$ for $N>0$ and $f(N)>g(N)$ for 
$0<N<K$, $f(N)<g(N)$ for $N>K$ ($K$ is the carrying capacity of the environment), the 
positive equilibrium $K$ is stable: all positive solutions converge to 
$K$ and are monotone.
It was argued that the observed data usually oscillates about the carrying 
capacity; in order to model this phenomenon, it was suggested to introduce
delay in the production term
\begin{equation}
\label{02}
\frac{dN}{dt}=f(N(t-\tau))-g(N);
\end{equation}
the latter equation can have oscillatory solutions, and the delay 
incorporated in the right hand side can be interpreted as maturation,
production or digestion effects. It is usually assumed that the 
mortality rate was proportional to the present population level 
\begin{equation}
\label{02new}
\frac{dN}{dt}=-\mu N(t)+ f(N(t-\tau)),~~\mu>0.
\end{equation}
The global behavior of solutions of (\ref{02new}) has been extensively 
studied in literature, in particular in the cases of negative and 
positive feedback
(see, for example, \cite{KWW1999, KW2001} and references therein), 
the chaotic behavior is impossible in the case of the monotone feedback
\cite{MPS}. However in the case when $f(x)$ is a unimodal
function, i.e., increases for $x<K$ and decreases for $x>K$, there may be 
delay induced instability and complex dynamics \cite{MG1,LW}. 
For a detailed overview of the literature on the dynamics of (\ref{02new}) 
see the recent papers \cite{RostWu,RostLiz}. 
It is demonstrated in \cite{RostWu,RostLiz} that if $f$ is a unimodal function
and positive equilibrium $K$ of the equation $x_{n+1}=f(x_n)$ is globally 
asymptotically stable, then all solutions of (\ref{02new}) tend to $K$. In 
particular, if $f$  has a negative
Schwarzian derivative, then local stability of the equilibrium of the
difference equation implies its global attractivity \cite{Singer}.
To the best of our knowledge, the first
delay-independent stability conditions were obtained in \cite{IvShark}.
In the present paper we will try to answer the general question: what are 
intrinsic properties of the reproduction function $f$ which allow us to 
conclude
that any solution of the equation with a finite memory converges to the
equilibrium?
Here we consider both general delays (including integral terms) and 
continuous functions $f$ which may have multiple extrema, tend to infinity 
at infinity etc.

As special cases, (\ref{02new}) includes 
the Nicholson's blowflies equation \cite{nichol,GBN} and the Mackey-Glass 
equation \cite{MG1,MG2}.
The Nicholson's blowflies equation
\beq{71}
\dot{x}(t)= -\delta x(t)+  px(t-\tau) e^{-a x(t-\tau)}
\eeq
was used in \cite{GBN} to describe the periodic
oscillation in Nicholson's classic experiments \cite{nichol} with the
Australian sheep blowfly, {\em Lucila cuprina}.
Equation \rf{71} with a distributed delay was studied in 
\cite{CAMQ2006}, where comprehensive results were obtained for the case $\delta < p < 
\delta e$.

The Mackey-Glass equation \cite{MG1,MG2}
\beq{72}
\dot{x}(t) = \frac{ax(t-\tau)}{1+x^{\gamma} (t-\tau)}-bx(t)
\eeq
models white blood cells production. Local and global stability of
the positive equilibrium for equation (\ref{72}) with variable delays
was studied in \cite{Liz1,Liz2,KubSak,Saker,BB2006,BBDCDIS}; to the best 
of our knowledge, there are no publications on  (\ref{72}) with a 
distributed delay.

To incorporate random environment influence, some 
authors included noise in (\ref{71}) and studied attractivity conditions, 
see, for example, \cite{Shaikhet}. 
However, in applied problems not only the derivative but also the delay 
value can be perturbed. We assume that the production delay is not a 
constant $\tau$ but some distributed value which leads to the equation
\beq{first}
\dot{x}(t) =  r \left[ \int\limits_{-\infty}^t f( x(s)) d_s R(t,s) - x(t) 
\right],
\eeq
where ${\displaystyle \int_{t-a}^{t-b} d_s R(t,s) }$ is the probability 
that at 
time $t$ the maturation delay in the production function is between $b$ 
and $a$, where $0<b<a$. We will assume that very large delays are 
improbable, substituting $-\infty$ in the lower bound with $h(t)\leq t$
which tends to infinity as $t \to \infty$.
In the present paper we consider a rather general form of $f$, which 
includes unimodal functions, as well as functions with several extrema.
The only requirement is that $f(x)$ has the only positive fixed point.
The main result claims that if this fixed point is a global attractor
for all positive solutions of the difference equation
\beq{first_difference}
x_{n+1}=f(x_n)
\eeq
then all solutions of \rf{first} with positive initial conditions tend to 
this fixed point as well.
To some extent this establishes a link between stable 
differential equation \rf{01}  and difference equation 
\rf{first_difference} which can undergo a series of bifurcations and even 
transition to chaos.  
If \rf{first_difference} is globally stable, so is \rf{first}. 
If the unique positive equilibrium of \rf{first_difference} is unstable,
\rf{first} can be stable or not, depending on the delay.

The paper is organized as follows.
In Section 2 we prove that all solutions with positive initial conditions 
are positive and bounded and establish some estimates for the lower
and the upper bounds. Section 3 presents sufficient conditions under which 
all positive solutions converge to the positive equilibrium.
In Section 4 delay-dependent stability is investigated. In particular, it 
is demonstrated that equations are globally attractive for delays small 
enough; if \rf{first_difference} is unstable, then we can find such delays 
that the positive equilibrium of \rf{first} is not a global attractor.
In Section 5 these results are applied to equations of population dynamics 
with a unimodal reproduction function and a distributed delay, in 
particular, to the Nicholson's blowflies and Mackey-Glass equations;
some open problems are presented.

\section{Boundedness and Estimates of Solutions}

We consider the equation with a distributed delay 
\beq{1a}
\dot{x}(t) =  r(t) \left[ \int_{h(t)}^t f( x(s)) d_s R(t,s) - x(t) 
\right], ~t \geq  0,
\eeq
and the initial condition
\beq{2star}
x(t)=\varphi(t), ~ t\leq 0.
\eeq
As special cases, \rf{1a} includes
\begin{enumerate}
\item
{\bf The integrodifferential equation} 
\beq{int1}
\dot{x}(t)= r(t) \left[ \int_{h(t)}^t K(t,s) f(x(s))\, ds - x(t) 
\right]
\eeq
corresponding to the absolutely continuous $R(t, \cdot)$ for any $t$.
Here
$$\int_{h(t)}^t K(t,s) ~ds=1 \mbox{~~for any~~~} t, ~~
K(t,s) = \frac{\partial}{\partial s} R(t,s) \geq 0 $$ is defined
almost everywhere.
\item
{\bf The equation with several concentrated delays}
\beq{conc1}
\dot{x}(t)= r(t) \left[ \sum_{k=1}^m a_k(t) f\left( x(h_k(t))
\right) - x(t) \right],
\eeq
with $a_k(t) \geq 0$, $k=1, \cdots , m$, where ${\displaystyle 
\sum_{k=1}^m
a_k(t) = 1}$ for any $t$.
This corresponds to ${\displaystyle R(t,s)=\sum_{k=1}^m a_k(t)  
\chi_{(h_k(t),\infty)} (s) }$, where $\chi_{I}(t)$ is the characteristic
function of interval $I$.
\end{enumerate}

\noindent
{\bf Definition.} An absolutely continuous in $[0, \infty)$
function $x: \RR \rightarrow \RR$ is
called {\em a solution of the problem}
(\ref{1a}),(\ref{2star}) if it satisfies equation (\ref{1a})
for almost all $t\in [0,\infty)$
and conditions (\ref{2star}) for $t\leq 0$.
\vspace{2mm}

The integral in the right hand side of \rf{1a}  
should exist almost everywhere. In particular, for \rf{int1}
with a locally integrable kernel, $\varphi$
can be any Lebesgue measurable essentially bounded function. For
\rf{conc1} $\varphi$ should be a Borel measurable bounded function. For
any distribution $R$ the integral exists if $\varphi$ is bounded and
continuous (here we assume $f$ is continuous). Besides, as is commonly
set in population dynamics models, $\varphi(t)$ is 
nonnegative and the value at the initial point is positive. 

Consider (\ref{1a}),(\ref{2star}) under the following assumptions.
\begin{description}
\item{{\bf (a1)}} 
$f:[0,\infty) \to [0,\infty)$ is a continuous function satisfying 
Lipschitz condition $|f(x)-f(y)| \leq L |x-y|$, $x,y \geq 0$, 
$f(0)=0$, $f(x)>x$ for $0<x<K$ and $0<f(x)<x$ for $x>K$;
\item{{\bf (a2)}} 
$h:[0,\infty)\rightarrow \RR$, is a Lebesgue measurable
function, $ h(t)\leq t$,
$\lim\limits_{t\rightarrow \infty}h(t)=\infty;$
\item{{\bf (a3)}} 
$r(t)$ is a Lebesgue measurable essentially bounded on $[0,\infty)$ 
function, $r(t) \geq 0$ for any $t \geq 0$, ${\displaystyle 
\int_0^{\infty} r(s)~ds = \infty}$;
\item{{\bf (a4)}} 
$R(t, \cdot)$ is a left continuous nondecreasing function
for any $t$, $R(\cdot,s)$ is locally integrable for
any $s$, $R(t,s)=0$, $s \leq h(t)$, $R(t,t^+)=1$.
Here $u(t^+)$ is the right side limit of function $u$ at point $t$.
\item{{\bf (a5)}} 
$\varphi: (-\infty,0] \to \RR$ is a continuous bounded function,
$\varphi(t) \geq 0$, $\varphi(0)>0$.
\end{description}

First, let us justify that the solution of \rf{1a},\rf{2star} exists and 
is unique.

Denote by ${\bf L}^2([t_0,t_1])$ the space of Lebesgue measurable 
functions $x(t)$ such that
$$Q=\int_{t_0}^{t_1} (x(t))^2~dt<\infty , ~~\| x 
\|_{{\bf L}^2([c,d])}=\sqrt{Q},$$
by ${\bf C}([t_0,t_1])$  the space of continuous in
$[t_0,t_1]$ functions with the $\sup$-norm.

We will use the following result from the book of Corduneanu
\cite{Cordun} (Theorem 4.5, p. 95). We recall that operator $N$ is
{\em causal} (or {\em Volterra}) if for any two functions $x$ and $y$
and each $t$ the fact that $x(s)=y(s)$, $s \leq t$, implies
$(Nx)(s)=(Ny)(s)$, $s \leq t$.

\begin{uess}
\label{lemma1} \cite{Cordun}
Consider the equation
\beq{55}
\dot{y}(t)=({\cal L}y)(t)+({\cal N}y)(t), ~~t\in [t_0,t_1],
\eeq
where ${\cal L}$ is a linear bounded causal operator,
$N: {\bf C}([t_0,t_1]) \to {\bf L}^2([t_0,t_1])$ is a nonlinear causal 
operator which
satisfies
\beq{56}
\| {\cal N}x - {\cal N}y \|_{{\bf L}^2([t_0,t_1])} \leq \lambda \| x-y 
\|_{{\bf 
C}([t_0,t_1])}
\eeq
for $\lambda$ sufficiently small. Then there exists a unique absolutely
continuous solution of \rf{55} in $[t_0,t_1]$, with the initial function 
being equal to zero for $t<t_0$.
\end{uess}

\begin{guess}
\label{theorem0}
Suppose (a1)-(a5) hold. Then there exists a unique solution of
(\ref{1a}),\rf{2star}.
\end{guess}
{\bf Proof.} To reduce \rf{1a} to the equation with the zero initial
function, for any $t_0 \geq 0$ we can present the integral as a sum of two
integrals 
\beq{57}   
\dot{x}(t)= -r(t) x(t) + r(t) \int_{t_0}^t f(x(s))~d_s R(t,s) +
r(t) \int_{t_0}^t f(\varphi(s))~d_s R(t,s),
\eeq
where
$$ x(t)=0,~t<t_0, ~~~\varphi(t)=0,~t \geq t_0.$$
Here $t_0 \geq 0$ is arbitrary, so we begin with $t_0=0$ and proceed to a
neighboring $t_1$ to prove the existence of a local solution.
Then in \rf{55}
$$({\cal L}x)(t)=-r(t) x(t), ~~({\cal N}x)(t)= r(t) \int_{t_0}^t 
f(x(s))~d_s R(t,s) + F(t),$$
where 
$$
F(t)=r(t) \int_{t_0}^t f(\varphi(s))~d_s R(t,s), ~~
 x(t)=0,~t<t_0, ~~~\varphi(t)=0,~t \geq t_0,$$
and for any $\lambda>0$ there is $t_1$, such that
\begin{eqnarray*}
\| {\cal N}x - {\cal N}y \|_{{\bf L}^2([t_0,t_1])} & \leq & |r(t)| \left\| 
\int_{t_0}^t
|f(x(s))-f(y(s))|~d_s R(t,s) \right\|_{{\bf L}^2([t_0,t_1])}  \\
& \leq & L \, {\rm ess}\sup_{t \geq 0} |r(t)| \max_{s \in [t_0,t_1]} 
|x(s)-y(s)| \left\|  \int_{t_0}^t d_s
R(t,s) \right\|_{{\bf L}^2([t_0,t_1])} \\
& \leq & L \, {\rm ess} \sup_{t \geq 0} |r(t)| \, \| x(s)-y(s) 
\|_{C([t_0,t_1])} 
|t_0 - t_1|
 \leq \lambda \| x-y \|_{{\bf C}([t_0,t_1])}
\end{eqnarray*}
for ${\displaystyle |t_0-t_1| \leq \lambda/(L \, {\rm ess} \sup_{t \geq 0} 
|r(t)| )}$, 
where $L$ was defined in (a1), here $\lambda$ can be chosen small enough.
By Lemma \ref{lemma1} this implies existence and uniqueness of a local 
solution for \rf{1a}. This solution is either global or there exists $t_2$ 
such that either
\beq{58}
\liminf_{t \to t_2} x(t) = -\infty
\eeq
or
\beq{59}
\limsup_{t \to t_2} x(t) =\infty ~.
\eeq
The initial value is positive, so as far as $x(t)>0$, the solution
is not less than the solution of the initial value problem 
$\dot{x}+r(t)x=0$, $x(0)=x_0>0$ which is positive and the former case 
\rf{58} is impossible.
In addition, $\dot{x}(t)<0$ for any $$ x(t)>K,~x(t) \geq \max\left\{ 
\max_{0 \leq s \leq t} x(s), \sup_{s \leq 0} \varphi(s) \right\}, $$ which 
contradicts
\rf{59}. Thus there exists a unique global solution, which completes 
the proof.
\qed

\begin{guess}
\label{theorem1}
Suppose (a1)-(a5) hold.
Then the solution of (\ref{1a}),\rf{2star} is
positive for $t \geq 0$.
\end{guess}
{\bf Proof.} After the substitution
\beq{subst}
y(t)=x(t)\exp\left\{ \int_0^t r(\zeta)~d\zeta \right\},
\eeq
equation (\ref{1a}) becomes
\beq{3}
\dot{y}(t) = r(t) \exp\left\{ \int_0^t r(s)~ds\right\} \int_{h(t)}^t f 
\left( y(s) \exp\left\{ -\int_0^s 
r(\zeta)~d\zeta \right\} \right)~d_s R(t,s),
~t \geq 0. 
\eeq
Thus $y(0)>0$ and $\dot{y}(t) \geq 0$ as far as $y(s) \geq 0$, $s \leq t$,
consequently, $y(t)>0$ for any $t \geq 0$. Since the signs of $y(t)$
and $x(t)$ coincide, then $x(t)>0$ for any $t \geq 0$.
\qed

\noindent
{\bf Definition.} The solution $x(t)$ of (\ref{1a}),\rf{2star} is {\em
permanent} if there exist $A$ and $B$, $B\geq A>0$,  such that
$$A \leq x(t) \leq B, ~t \geq 0.$$

In the following we prove permanence of all solutions  
of (\ref{1a}) with positive initial conditions; moreover, we  
establish bounds for solutions.

\begin{guess}
\label{theorem1a}
Suppose (a1)-(a5) hold. Then a solution of \rf{1a},\rf{2star} 
is permanent.
\end{guess}
{\bf Proof.} By Theorem \ref{theorem1} the solution is positive for $t
\geq 0$. By (a2) there exists $t_0>0$ such that $h(t)>0$, $t \geq t_0$.
Since the solution is a continuous positive function, then we can define
\beq{defxminmax}
x_{\min} = \min_{t \in [0,t_0]} x(t)>0, ~~
x_{\max} = \max_{t \in [0,t_0]} x(t).
\eeq
Without loss of generality we assume $x_{\min}<K$, $x_{\max}>K$;
otherwise, we can choose $\min\{ x_{\min},\lambda K \}$, $\max\{ x_{\max}, 
K/\lambda\}$, where $0<\lambda<1$, as $x_{\min}$ and $x_{\max}$, 
respectively.
By (a1) the following values are positive
\beq{deffminmax}
M = \max_{x\in [x_{\min},x_{\max}]} f(x), ~~
m = \min_{x\in [x_{\min},x_{\max}]} f(x).
\eeq
Define 
\beq{defAB}
B= \max \left\{ M, x_{\max}, \max_{x \in [0,K]} f(x) \right\},
%\eeq
%\beq{defB}  
~~~
A= \min \left\{ m, x_{\min}, \min_{x \in [K,B]} f(x) \right\}.
\eeq
Since $f(x)>x$, $0<x \leq A$ and $f(x)<x$, $x \geq B$, then 
there exists $\delta>0$ such that
$f(x) \geq A$ for $A-\delta \leq x \leq B$ and $f(x) \leq B$ for $A \leq x 
\leq B+ \delta$.
Let us demonstrate
\beq{bounds}
x(t) \in [A,B], ~t \geq 0.
\eeq
By the definition of $A,B$ we have $x(t_0) \in [A,B]$. 
Suppose the contrary: $x(t)>B$ or $x(t)<A$ for some 
$t>t_0$.

First, let $x(t)>B$ for some $t>t_0$. Then $x(t)=B+\varepsilon$ for some 
$\varepsilon \leq \delta$. 
Denote
$$S_1= \{ t >t_0 | x(t)>B+\varepsilon\}, ~~t^{\ast} = \inf S_1.$$
Since $x(t_0) \leq B$ then the set 
$$ S_2= \{ t | t_0 \leq t < t^{\ast}, ~x(t) \leq B\}
$$
is nonempty, denote $t_{\ast} = \sup S_2$. Then $x(t_{\ast})=B$, 
$x(t^{\ast})=B+\varepsilon$; we also have $B\leq x(t) \leq B+\delta$ and $f(x) \leq 
B$ in the interval 
$[t_{\ast},t^{\ast}]$, thus the derivative is nonpositive
$$
\dot{x}(t)  =  r(t) \left[ \int_{h(t)}^t f(x(s))~d_s R(t,s) -
x(t) \right] \leq r(t) \left[ \int_{h(t)}^t B~d_s R(t,s) - B
\right]=0,
$$
which contradicts the assumption $x(t^{\ast})=B+\varepsilon > B= x(t_{\ast})$.

Similarly, let us assume that $x(t)=A-\varepsilon$ for some 
$\varepsilon>0$, $\varepsilon<\delta$ and some $t>t_0$. After introducing
$$S_1= \{ t >t_0 | x(t)<A-\varepsilon\}, ~t^{\ast} = \inf S_1,
~~S_2= \{ t | t_0 \leq t < t^{\ast}, ~x(t) \geq A \}, ~t_{\ast} = \sup 
S_2,
$$
we have $x(t_{\ast})=A$,
$x(t^{\ast})=A-\delta$, $A-\delta\leq x(t) \leq A$ and $f(x(t)) \geq A$ 
for
$t \in [t_{\ast},t^{\ast}]$, hence
$$  
\dot{x}(t)  =  r(t) \left[ \int_{h(t)}^t f(x(s))~d_s R(t,s) -
x(t) \right] \geq r(t) \left[ \int_{h(t)}^t A~d_s R(t,s) - A
\right]=0.
$$
This contradicts the assumption $x(t^{\ast}) < x(t_{\ast})$.
Consequently (\ref{bounds}) is valid for $t\geq t_0$ and thus for any 
$t\geq 0$, the bounds are positive, so the solution is permanent, which 
completes the proof.
\qed

\noindent
{\bf Example 1.} The statement of Theorem \ref{theorem1a} is not valid if
we omit the condition ${\displaystyle \lim_{t \to \infty} h(t)=\infty}$.
Consider the equation
\beq{ex1eq1} 
\dot{x}(t) = 5x(h(t)) e^{-x(h(t))}-x(t),~t\geq 0, ~~h(t) \equiv -1, 
~~x(t)=t+1,
~t\in [-1,0],
\eeq
which is equivalent to the initial value problem
\beq{ex1eq2}
\dot{x}(t)+x(t)=0, ~~x(0)=1,
\eeq
its solution ${\displaystyle x(t)=e^{-t}}$ tends to zero as $t \to \infty$ 
and so is not permanent.

\section{Absolute Global stability for Stable Difference Equations}

One of the main steps in establishing global stability property is the 
proof of the fact that all nonoscillatory about the equilibrium solutions 
tend to this equilibrium (see, for example, \cite{GL}). For 
ordinary differential equations all solutions are nonoscillatory, for 
retarded equations it depends on the delay. Below we demonstrate that 
convergence of nonoscillatory solutions to the equilibrium is quite a 
common property which is valid for any reproduction function with a unique
positive equilibrium. It can be interpreted as: ``if nonoscillatory, 
solutions of delay equations behave asymptotically similar to ordinary
differential equations".   
\vspace{1mm}

\noindent
{\bf Definition.} A solution $x(t)$ of (\ref{1a}),\rf{2star} is {\em 
nonoscillatory about $K$} if there exists $\tau>0$ such that 
either $x(t)>K$ or $x(t)<K$ for all $t \geq \tau$. Otherwise, $x(t)$ {\em
oscillates about $K$}.

\begin{guess}
\label{nonosciltheor}
Suppose (a1)-(a5) hold.
Any nonoscillatory about $K$ solution of (\ref{1a}),\rf{2star} converges 
to $K$.
\end{guess}
{\bf Proof.}
First, let $x(t)<K$, $t \geq \tau$. Without loss of generality we can 
assume $\tau=0$. By (a2) there exists $t_0 \geq 0$ such that $h(t) \geq 0$ 
for $t \geq t_0$. Denote $A$ as in \rf{defAB}.
By Theorem \ref{theorem1a} we obtain that $x(t) \geq 
A$ for any $t \geq 0$. Since $f$ is continuous and $f(x)>x$ for $x<K$, then
${\displaystyle m_0=\inf_{x \in [A,K]} f(x) > A }$. 
%Further, assuming 
%$x(t)=A_0<A$ for some $t>t_0$ we denote by $t^{\ast}$ the first $t$ where
%$x(t)=A_0$ and by $t_{\ast}$ the largest $t<t^{\ast}$ where $x(t)=A$.
%Then $x(t)$ increases in $[t_{\ast},t^{\ast}]$, while
%$A=x(t_{\ast})>x(t^{\ast})=A_0$ which is a contradiction.
There may be two 
cases: $m_0=K$ and $m_0<K$. In the former case, since $x(t)<K$, we have
$$\dot{x}(t)>r(t)\left[ \inf_{x \in [A,K]} f(x) - x(t) \right] > 0,
~t \geq t_0,
$$
as far as $x(t)>A$, thus the solution of the delay differential 
equation is not less than the solution of $\dot{x}(t)=r(t)[K-x(t)]$, 
$0<x(t_0)<K$, which is increasing and by (a3) (the integral of $r(t)$ 
diverges) tends to $K$.

\begin{figure}[ht]
\centering
\includegraphics[scale=0.7]{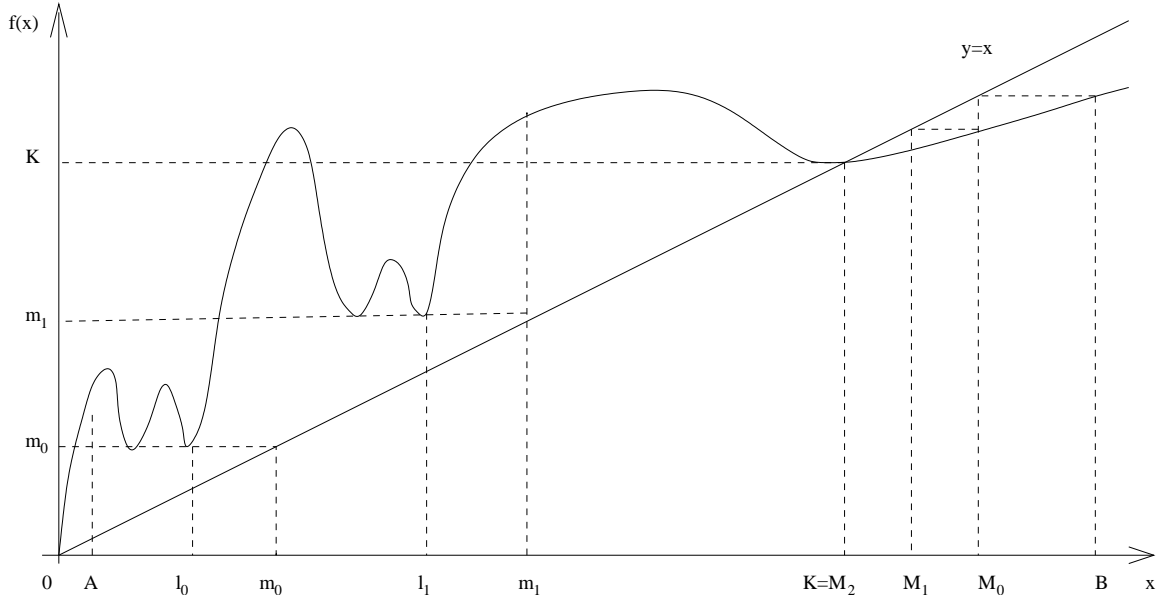}
\caption{For an arbitrary reproduction function with one positive
equilibrium $K$ we construct a series of such points that eventually
a solution is in $[m_j, K]$, if it does not exceed $K$ and
is in $[ K, M_j]$ if a solution is not less than $K$.
}
\label{figure1}   
\end{figure}

Consider the latter case $m_0<K$. 
%Since $A>0$ (the solution is 
%persistent, i.e., the lower bound exceeds zero), then $m_0>0$. Similar to 
%the previous argument, we obtain 
%that $x(t) \geq A$ for any $t \geq 0$.
By the definition of $m_0$ and $f(x)>x$ for $x \in [A,K]$ we have 
${\displaystyle l_0=\sup\{x<K| f(x)\leq m_0\} < m_0 }$. Taking any
$\alpha$, $l_0<\alpha<m_0$ and assuming $x(t) \leq \alpha$ for any $t$, we 
obtain
$$ \dot{x}(t)= r(t) \left[ \int_{h(t)}^t f(x(s)) ~d_s R(t,s) - x(t)
\right] \geq r(t)(m_0-\alpha)>0,$$
which leads to a contradiction $x(t)
\to \infty$ as $t\to \infty$ since $\int_0^{\infty} r(s)~ds $ diverges.
Thus, $x(t_*) \geq \alpha$ for some $t_*$; 
moreover, since $\dot{x}(t) \geq 0$ as $x(t) \leq m_0$ 
then $x(t) \geq \alpha$ for any $t \geq t_*$.
Let $h(t) \geq t_*$, $t>t^*$ for some $t^*$.
%Thus, $x(t)\geq \alpha$ for some $t_{\ast}$; moreover, since $ 
%\dot{x}(t)\geq 0$ as
%$x(t) \leq m_0$ then $x(t)\geq l_1$ for any $t \geq t_{\ast}$. 
%Let $h(t)\geq \alpha$ for $t \geq t^{\ast}$. 
By the definition of $l_0$ and $\alpha$ we have ${\displaystyle \tilde{m}=
\inf\limits_{x\in [\alpha,K]} f(x)>m_0}$ and
as far as  $x(t) \leq m_0$ and $h(t) \geq t^{\ast}$ the following 
inequality holds
$$ \dot{x}(t)= r(t) \left[ \int_{h(t)}^t f(x(s)) ~d_s R(t,s) - x(t) 
\right] \geq r(t)(\tilde{m}-m_0)>0.$$
Assuming $x(t) \leq m_0$ for any $t$ we again obtain a contradiction.
Thus, there exist $\mu_1$ and $t_1 > t_0$ such that $x(\mu_1) \geq m_0$ and 
$h(t) \geq \mu_1$ for $t \geq t_1$. Then $x(t) \geq 
m_0$ for any $t \geq \mu_1$ and $x(h(t))\geq m_0$, $t \geq t_1$.

Further, let ${\displaystyle m_1 = \inf_{x \in 
[m_0,K]} f(x)<K}$, here $m_0<m_1$ since $f(x)>x$ for $0<x<K$. Similarly, 
we obtain $x(t)\geq m_1$ whenever $t>t_2$, for some $t_2>t_1$. 

We continue this 
process. It can be finite (for example, in Fig. \ref{figure1} we have $m_2=K$, where
the process stops and we deduce $x(t) \to K$ as $t \to \infty$) or
infinite (see the branch $x(t)>K$ of Fig. 1). In the infinite case we have 
an increasing sequence $\{ m_j \}$, ${\displaystyle m_{j+1}=\min_{x 
\in [m_j,K]} f(x) }$ which does not exceed $K$, so this sequence has a 
limit $d$. Since $f(x)$ is continuous then ${\displaystyle 
d=\min_{x\in [d,K]} f(x) }$. If $d<K$ then $f(x)-d$ should attain its 
minimum in $[d,K]$ but $f(x)>x \geq d$, so this minimum is positive and
the equality ${\displaystyle d=\min_{x\in [d,K]} f(x) }$ leads to a 
contradiction.

Further, let $x(t)>K$. Similarly, we define $B$ as in \rf{defAB}
and
${\displaystyle M_0=\max_{x \in [K,B]} f(x)<B }$. There may be two   
cases: $M_0=K$ and $M_0>K$. In the former case we obtain $x(t) \to K$ as
$t \to \infty$.
Consider the latter case. By Theorem \ref{theorem1a}
we have $x(t) \leq B$ for any $x \geq 0$.
By the definition of $M_0$ and $f(x)<x$, $x \in [K,B]$ we have
${\displaystyle s_0=\inf\{x>K| f(x)\geq M_0\} > M_0 }$. 
Similar to the case $x(t)<K$ we demonstrate that there exists $\nu_1$ such 
that $x(\nu_1) \leq M_0$ and $x(t) \leq M_0$ for any $t \geq \nu_1$. Let 
$t_1$ be such that $h(t) \geq \nu_1$ for $t \geq t_1$. Further, we define
${\displaystyle M_1=\max_{x \in [K,M_0]} f(x)<M_0 }$.
We continue this process, it can be finite or infinite (Fig.~\ref{figure1}
illustrates an infinite process for $x(t)>K$).
Similar to the case $x(t)<K$ we obtain $x(t) \to K$ as $t \to \infty$, 
which completes the proof.
\qed

\noindent
{\bf Example 2.} Let us note that in the case of infinite delays 
nonoscillatory solutions do not necessarily tend to the positive 
equilibrium. For example, the solution of the equation
$$\dot{x}(t) + x(t) = \frac{e}{2}\, x(0) \, e^{-x(0)}, ~~x(0)=2,$$
which is ${\displaystyle x(t)= \frac{1}{e} + \left( 2-\frac{1}{e} \right) 
e^{-t}}$, tends to $1/e$ while the positive equilibrium is $1- \ln 2$,
the monotone solution is nonoscillatory.
\vspace{2mm}

Thus for any reproduction function with a unique fixed point $f(x)=x$ 
nonoscillatory solutions tend to the equilibrium; this is not, generally,
true for oscillatory solutions.
\vspace{2mm}

\noindent
{\bf Example 3.} Consider the Nicholson's blowflies equation
\beq{ex3eq1}
\dot{x}(t)= -\delta x(t)+  px(t-\tau) e^{-a x(t-\tau)}.
\eeq
Denote
\beq{ex3eq2}
\tau_k= \frac{1}{\delta \sqrt{aN^{\ast} (aN^{\ast}-2)}}
\left[ \arcsin 
\left( \frac{\sqrt{aN^{\ast}(aN^{\ast}-2)}}{aN^{\ast}-1} \right) + 2\pi k
\right], ~~k=0,1,2, \cdots ~,
\eeq
where $N^{\ast}=1/a \ln(p/\delta)$ is a positive equilibrium.
If $p>\delta e^2$ then 
the positive equilibrium is locally asymptotically stable
for $\tau \in [0, \tau_0)$ and is unstable (locally, thus it cannot be
globally attractive) for $\tau>\tau_0$, and \rf{ex3eq1} undergoes a Hopf 
bifurcation at $N^{\ast}$ when $\tau=\tau_k$ \cite{Li} for $~k=0,1,2, 
\cdots$.
\vspace{2mm}

Now we prove that absolute (delay-independent) convergence holds in some 
special cases. 

\begin{uess}
\label{lemmaimp}
Suppose (a1)-(a5) and at least one of the following conditions holds:
\begin{enumerate}
\item
$f(x)<K$ for any $0<x<K$.
\item
Denote by $x_{\max}([a,K])$, $0 \leq a \leq K$, the greatest point in $[a,K]$ 
where ${\displaystyle \max_{x \in [a,K]} f(x)}$ is attained and assume
\begin{equation}
\label{attractivitycond}
f(x_{\max}([a,K])) > K \mbox{~for some ~} a \in (0,K) 
\Rightarrow \!\!
\min_{ x \in [K, f(x_{\max}([a,K]))]} \!\!\!\! f(x) >a.
\end{equation} 
\end{enumerate}
Then any solution of (\ref{1a}),\rf{2star} converges
to $K$.
\end{uess}
{\bf Proof.} 
1. By Theorem \ref{theorem1a} any solution is permanent: $A \leq x(t) \leq 
B$, where $A<K<B$, $t \geq 
0$, and for some $t_0 \geq 0$ we have $h(t) \geq 0$, $t\geq t_0$. Denote 
${\displaystyle M_0=\max_{x \in 
[K,B]}f(x)}$ and ${\displaystyle m_0=\min_{x \in [A,K]}f(x)}$. 

If $M_0 \leq K$ then the derivative is negative for any $x(t)>K$ and
the solution either eventually does not exceed $K$ or is 
decreasing for any $t \geq t_0$. In the latter case the solution tends to
the equilibrium; if it has a different limit, we obtain that the 
derivative is less than a negative number, which is a contradiction.
In the former case, if there exists $t>0$ such that $x(t)<K$ (otherwise,
we have a nonoscillatory case $x(t) \geq K$ where we have already proved 
convergence), then $x(t) \leq K$ for any $t$ (assuming $x(t^{\ast})>K$ we 
obtain that the derivative of $x(t)$ is negative almost everywhere, while 
the function changes from $K$ to $x(t^{\ast})>K$), which is again a 
nonoscillatory case, by Theorem \ref{nonosciltheor} solution $x(t)$ 
converges to $K$. 

Thus, we can consider $M_0 >K$ only. Then we introduce $$ M_1=\max_{x \in 
[K,M_0]} f(x)<M_0, ~~ m_1 =\min_{x \in [m_0,K]}>m_0, ~~s_1=\min_{x \in 
[K,M_0]} f(x)
$$
and similarly define $M_j,m_j$ and $s_j$, $j=2,3,\cdots $.
There exists $\tau_1$ such that $\min\{ m_1,s_1 \} \leq x(t) \leq
M_1$ for $t\geq \tau_1$ and $t_1$ such that $h(t) \geq \tau_1$, $t\geq 
t_1$. Similar to the proof of Theorem \ref{nonosciltheor}  we obtain that 
there exists a sequence $\tau_1 \leq \tau_2 \leq \cdots \leq \tau_j \leq 
\cdots $ such that 
$$\min\{ m_j,s_j \} \leq x(t) \leq M_j \mbox{~~for~~~} t \geq \tau_j.$$
We assume that all $M_j>K$, otherwise we have an eventually monotone case.
Since for any continuous $f$ all three sequences are monotone 
(nonincreasing $M_j$ and nondecreasing $m_j$, $s_j$) and tend to $K$ (see 
the end of the proof of  Theorem \ref{nonosciltheor}), then 
%for any $\varepsilon>0$ there exists
%$\tau_l$ such that
%$$ K-\varepsilon \leq  x(t) \leq K+\varepsilon, ~~~t \geq \tau_l.$$
%Hence 
${\displaystyle \lim_{t \to \infty} x(t)=K}$.

2. Now suppose that \rf{attractivitycond} holds for any $a \geq 0$.
%Consider $x_0=x_{\max} [0,K]$. If $x_0=K$ then we just  obtain the case
%of item 1. Let us assume $x_0<K$ and demonstrate that there exists 
%$\tau_1$ such that $x_0 \leq x(t) \leq f(x_0)$, $t \geq \tau_1$.
We begin with $A \leq x(t) \leq B$ in Theorem \ref{theorem1a},  
%$A < x_0 <K< f(x_0) < B$, 
$t \geq 0$, where $A,B$ are defined in
\rf{defAB}, $h(t)>0$, 
$t \geq t_0$. Denote
\beq{addast1}
M_0=\max_{x \in [A,K]}f(x), ~ m_0=\min_{x \in [K,M_0]}f(x).
\eeq
The case $M_0 =K$ was considered in Part 1, thus we can restrict ourselves 
to the case $M_0>K$, maximum is attained at $$x_0=x_{\max}([A,K]),$$
here $x_0$ is the greatest point where the maximum is attained. 
Let us demonstrate that there exists $\mu_0$ such that $x_0 \leq x(t) 
\leq M_0$, $t \geq \mu_0$. 

If ${\displaystyle x(t) > S_0=\max \{ M_0, \max_{x \in [K,B]} f(x) \}=
\max_{x \in [A,B]} f(x) \geq K }$
for all $t$ then we have a nonoscillatory case and convergence to $K$,
which is a contradiction. Thus, $x(t_{\ast}) \leq S_0$ for some 
$t_{\ast}$; assuming there is $t^{\ast}>t_{\ast}$ such that $x(t_{\ast}) 
>S_0$ we obtain that the value of the function at the  end of the 
segment is higher than at the beginning point, while the derivative is 
nonpositive. So there exists $\tau_0 \geq t_0$ such that $x(t) 
\leq S_0$ for $t>\tau_0$. If $S_0 \neq M_0$ then ${\displaystyle S_1=
\max_{x \in [K,S_0]} f(x) < S_0}$. Similarly, we find $\tau_1$ such that
$x(t)\leq S_1$ for $t$ large enough. Since the sequence 
${\displaystyle S_{n+1}= \max_{x \in [K,S_n]} f(x)  }$ is nonincreasing 
and tends to $K$ then there exists $\tau_1$ such that $x(t) \leq M_0$,
$t \geq \tau_1$. Further, we will consider $t \geq t_1$, where
$h(t) \geq \tau_1$ whenever $t>t_1$, only.

If ${\displaystyle x(t) < s_0=\min \{ m_0, \min_{x \in [A,K]} f(x) \} =
\min_{x \in [A,M_0]} f(x)  }$
for all $t$, where $m_0$ was defined in \rf{addast1}, then we have a 
nonoscillatory case and convergence to $K$,   
which is a contradiction. As above, we prove that $x(t)\geq s_0$ for $t
\geq \nu_0$, here $\nu_0 \geq t_0$. If $s_0<m_0$, then we construct a 
sequence ${\displaystyle s_{n+1}= \min_{x \in [s_n,K]} f(x)  }$ which 
tends to $K$. If $m_0<K$ then some $s_j \geq m_0$. 
There is $\mu_0 \geq t_1$ such that $x(t)\geq m_0$, $t\geq 
\mu_0$.  Consequently, we have found $\mu_0$ such that $m_0 \leq x(t)  
\leq M_0$, $t \geq \mu_0$. Let $h(t) \geq \mu_0$, $t\geq t_2$.

Now, if $m_0=K$ then for $t\geq t_2$ the solution increases as far as
$x(t)<K$. Thus, either $x(t) < K$ for any $t\geq t_2$ and this monotone 
solution converges to $K$, or, if $x(t_{\ast}) \geq K$ for some $t_{\ast}$,
$x(t) \geq K$, $t \geq t_{\ast}$, and again we have a nonoscillatory
solution which converges to $K$.

Denote
$$M_1= \max_{x \in [m_0,M_0]}f(x), ~~m_1=\min_{x \in [m_0,M_1]}f(x).$$
Let us assume $M_1>K$, $m_1<K$ and demonstrate that there is $\mu_1 \geq t_2$ 
such that
$m_1 \leq x(t) \leq M_1$, 
$x(t) \geq \mu_1$, where
$M_1 = f\left( x_{\max}([m_0,K]) \right)$.
In fact, for any $t \geq t_2$ the solution is nonincreasing as far as
$x(t) \geq M_1$, which gives an upper bound.  Considering $t$ where the 
equation refers only to the values where this bound is valid, we obtain 
that the solution is nondecreasing if $x(t) \leq m_1$, 
which together with $m_1<K<M_1$ confirms the statement. We continue the 
induction process
$$
M_{n+1}=\max_{x \in [m_n,M_n]} f(x),~~
m_{n+1}=\min_{x \in [m_n, M_{n+1}]} f(x),  
$$
where $K \geq m_{n+1} > m_n $, $K \leq M_{n+1} < M_n$. This process can
be infinitely continued if all $m_n < K$, $M_n>K$ (otherwise, at certain 
stage we have a ``monotone" case which implies convergence), and there 
exists $\mu_j$, $j=0,1,2, \cdots$ such that $m_j \leq x(t) \leq M_j$, $t 
\geq \mu_j$. Let us assume that there are infinite sequences $\{ m_n \}$ 
and $\{ M_n \}$, both are monotone and bounded.
Then there exist limits
$$ x=\lim_{j \to \infty} m_j \leq K \leq \lim_{j \to \infty} M_j =X.$$ 
Since by the assumptions of the theorem $f(M_j)>m_j$, then $f(X) \geq x$.
If $x=f(X)$ then either $x=X=f(X)=X$ and $x(t)$ converges to $K$ or 
$\min\limits_{y \in [K,f_{\max}([x,K])]} f(y) \leq x$ which contradicts 
the assumptions of the theorem. 

Denote by $x_{\max,j}$, $x_{\min,j}$ the 
sequences where minimum $m_j$ and maximum $M_j$ are attained, respectively
(we recall that we choose maximal $x_{\max,j}$ and minimal $x_{\min,j}$
if an extremum is attained at several points). Then $f(x_{\max,j})=M_j$,
$f(x_{\min,j})=m_j$, $m_j \leq x_{\max,j} \leq K$, $K  \leq x_{\min,j}
\leq m_j$, $m_{j+1}>x_{\max,j}$, $M_{j+1}<x_{\min,j}$, which implies that
$x_{\max,j}$ and $x_{\min,j}$ tend to $x$ and $X$, respectively.
Let us prove that $x=f(X)$. 
Assume the contrary: $f(X)> x+ 2\varepsilon>0$.
Since ${\displaystyle \lim_{j \to \infty} x_{\min,j}=X}$ and $f$ is
continuous then there exist $n_0$ such that $f(x_{\min,j})>x+\varepsilon$,
$j \geq n_0$, but $f(x_{\min,j})=m_{j+1}$, so $m_{j+1}>x+\varepsilon$, 
which contradicts the equality ${\displaystyle \lim_{j \to \infty} m_j = 
\sup_j m_j=x}$. Thus ${\displaystyle \lim_{t \to \infty} x(t)=K}$,
which completes the proof.
\qed

Now we can prove the main result of the present paper.

\begin{guess}
\label{maintheorem1}
Suppose (a1)-(a5) hold and all positive solutions of the difference 
equation
\begin{equation}
\label{differ}
x_{n+1}=f(x_n)
\end{equation}
tend to $K$. 
Then any positive solution of (\ref{1a}),\rf{2star} converges
to $K$.
\end{guess}
{\bf Proof.} If either conditions of Part 1 of Lemma \ref{lemmaimp} or
\rf{attractivitycond} hold then $K$ attracts all positive solutions of 
(\ref{1a}),\rf{2star}. Let us assume that \rf{attractivitycond} does not 
hold (and we do not have a monotone case as in Part 1 of Lemma 
\ref{lemmaimp}), which means that
for some $a \in (0,K)$ we have 
$$
\min_{ x \in [K, f(a)]} f(x) \leq a, \mbox{~~or~~~} K < b  \leq f(a), ~~ 
f(b) \leq a, \mbox{~~where~~} b \in [K, f(a)].
$$
We can assume $f(b)=a$, otherwise, since $f(b)<a$, $f(K)=K>a$ 
and $f$ is continuous, then there is $c \in [K,b]$ such that
$f(c)=a$.  Thus, $b>K$ and $f^2(b)=f(f(b))=f(a) \geq b$. Since $f(x)<x$ 
for any 
$x>K$ then for any $x> M:=\max\limits_{x \in [0,K]} f(x)$ we have 
$f^2(x)<x$, 
then there is a fixed point of $f^2$ in the segment $[b,M]$,
in addition to a fixed point $x=K$. Since the fixed point $K$ 
of \rf{differ} cannot be a global attractor unless $K$ is the only fixed 
point of $f^2$ \cite{Coppel} then $K$ is not a global 
attractor of (\ref{differ}), which completes the proof. 
\qed

Here we have not considered the case when $f$ has no positive equilibria;
however for completeness we will consider this case as well.

\begin{guess}
\label{maintheorem2}
Suppose (a2)-(a5) hold, $f(0)=0$ and $f(x)<x$ for any $x>0$.  
Then any positive solution of (\ref{1a}),\rf{2star} converges
to zero.
\end{guess} 
{\bf Proof.} Let $t_0$ be such that $h(t) \geq 0$ for $t \geq t_0$.
Denote ${\displaystyle M_0= \sup_{t \in [0,t_0]} x(t) }$.
Since $f(x)<x$ then $x(t) \leq M_0$ for any $t \geq 0$.
The solution is decreasing as far as $x>M_1$, where ${\displaystyle
M_1=\max_{x \in [0,M_0]} f(x)<M_0}$.

\begin{figure}[ht]
\centering
\includegraphics[scale=0.95]{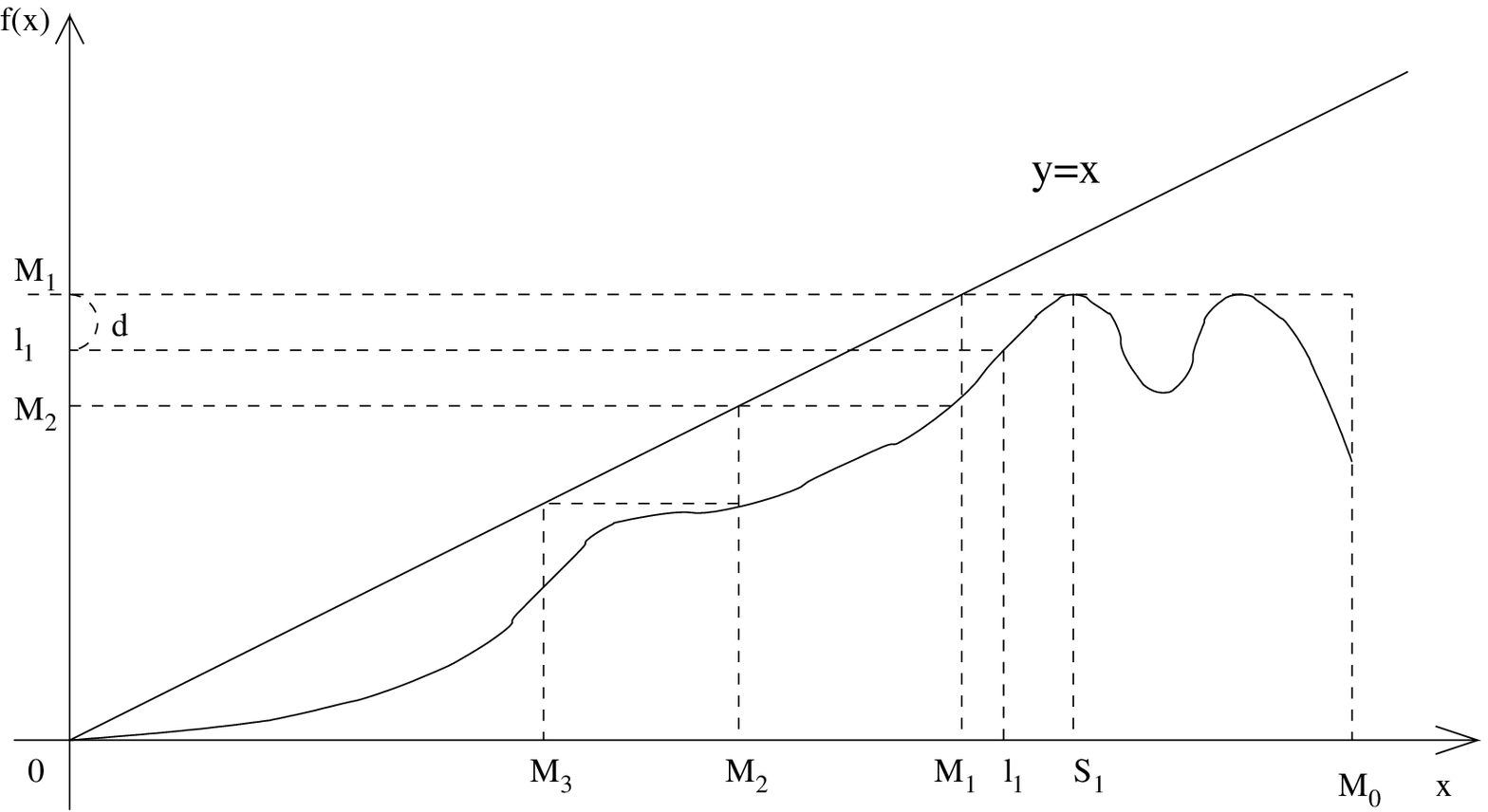}
\caption{If $f(x)<x$ for any $x>0$ then all positive solutions tend to zero.
}
\label{figure2}   
\end{figure}

Suppose $S_1$ is the smallest point not exceeding $M_0$ where this maximum  
is attained (see Fig.~\ref{figure2}). Since $f(x)<x$, then $S_1>M_1$. Let us 
choose $l_1$ such that
$M_1<l_1<S_1$; as far as $x(t) \geq l_1$ we have $\dot{x}(t) \leq
r(t)(M_1-l_1)$, thus eventually any solution is less than $l_1$. Further,
for some $d>0$ we have $M_1-f(x) \geq d$ if $0<x\leq l_1$. Thus the
derivative at any point $x(t)>M_1$ does not exceed (we recall that 
$x(t)\leq l_1$) the negative value of $-d\, r(t)$. Consequently, there 
exists $t_1$ such that $x(t) \leq M_1$ for $t \geq t_1$. 
By induction, we denote ${\displaystyle M_{k+1}=\sup_{x \in [0,M_k]} 
f(x)<M_k }$, $k=1,2, 
\cdots $, and prove that there exists $t_{k+1} >t_k$ such that $x(t) \leq 
M_{k+1}$ for $t>t_{k+1}$. Since ${\displaystyle \lim_{n \to \infty} M_n 
=0}$ then
${\displaystyle \lim_{t\to \infty} x(t)=0}$, i.e., any positive solution
tends to zero.
\qed

\section{The Case When Stability Is Delay-Dependent}

In this section we consider the case when stability properties 
of \rf{1a} depend on the delay. First, we prove that for small delays 
all solutions tend to positive equilibrium $K$.
Second, we demonstrate that if \rf{attractivitycond} does not hold then 
there exists equation \rf{1a} with parameters satisfying (a1)-(a5) such 
that its attracting set is as close to $[m,M]$ as prescribed, where
\beq{mM}
M=\max_{x\in[0,K]} f(x) > K, ~~m=\min_{x\in[K,M]} f(x)
\eeq 
and ${\displaystyle \max_{x\in[m,K]} f(x)=M}$ as well.

Let us note that the Lipschitz condition in (a1) implies
$| f(x) -K| \leq L |x-K|$.

\begin{guess}
\label{theorem5}
Suppose (a1)-(a5) hold and 
%there exists $\lambda$, $0<\lambda<1$ such that 
\begin{equation}
\label{cup}
\limsup_{t \to \infty} \int_{h(t)}^t r(s)~ds < \frac{1}{L+1}~.
\end{equation}
Then any positive solution of (\ref{1a}),\rf{2star} converges
to $K$.
\end{guess}
{\bf Proof.}
Without loss of generality we can assume that \rf{cup} is satisfied 
for $t \geq 0$. We can find a positive $\lambda<1$ such that 
${\displaystyle \int_{h(t)}^t 
r(s)~ds < \frac{\lambda}{L+1} }$ for $t \geq 0$. By Theorem
\ref{theorem1a} any solution is permanent: 
$$ K-a_0 \leq x(t) \leq K+b_0, ~~t \geq 0, ~0<a_0<K, ~b_0>0.$$
Let us denote $\alpha=\max \{ a_0, b_0 \}$ and prove that there exists 
$t_1 \geq 0$ such that $x(t) \leq K+\lambda \alpha$ for any $t > t_1$.
If solution $x(t)$ is nonoscillatory about $K$ then by Theorem 
\ref{nonosciltheor} it tends to the equilibrium and thus for 
$\varepsilon=\lambda \alpha>0$ there exists $t_1$ such that $x(t) \leq 
K+\varepsilon$ for $t \geq t_1$. Now assume that $x(t)$ is oscillating,
$t_0$ is such a point that $h(t) \geq 0$ for $t \geq t_0$ and 
$x(t_{\ast})=K$, $x(\tau_{\ast})=K$, where $x(t)>K$ for $t \in 
(t_{\ast},\tau_{\ast})$. A continuous function $x(t)$ attains its maximum
in $[t_{\ast},\tau_{\ast}]$, let ${\displaystyle x(t^{\ast})
= \max\limits_{s \in [t_{\ast},\tau_{\ast}]} x(s) = M_0 }$. We claim that 
$x(t^{\ast}) \leq K+ \lambda \alpha$. Assume the contrary: $x(t^{\ast}) 
>K+ \lambda \alpha$. 
By (a1) we have ${\displaystyle \max_{x \in [K,M_0]} f(x)=M_1<M_0}$. 
Denote  
$$M_2= \max \{ M_1, K+ \lambda \alpha \}, ~~
\tau = \sup \left\{ t \in [t_{\ast},t^{\ast}] | x(t) \leq M_2 
\right\}.$$
Then $x(t)> M_2 \geq K+\lambda \alpha$, $t \in [\tau,t^{\ast}]$.
 
%We will demonstrate 
%that also $h(t)>t_{\ast}$ for $t\in [\tau,t^{\ast}]$.
Since $K-\alpha \leq x(t) \leq K+\alpha$ for any $t \geq 0$ then by (a1)
$$
\dot{x} (t) = r(t) \left[ \int_{h(t)}^t f(x(s)) ~d_s R(t,s) - x(t) 
\right] \leq r(t) \left[  \max\limits_{x \in [K-\alpha, K+\alpha], x>0}
f(x) - (K-\alpha) \right] 
$$
\beq{star31}
= r(t) \left[ \max\limits_{x \in [K-\alpha, K+\alpha], x>0}
[f(x)-K] +\alpha \right] \leq r(t) [L\alpha+\alpha] = \alpha(L+1)r(t).
\eeq
We remark that $\dot{x} (t) <0$ whenever the following conditions hold:
$t \in (\tau, t^{\ast})$, $h(t) \geq t_{\ast}$ (all referred 
prehistory of the solution is between $K$ and $M_1<M_0$) and $x(t)>M_2 
\geq M_1$. Since $f(\tau)=M_2$, $f(t^{\ast})>M_2$, then there are points in  
$(\tau,t^{\ast})$ where the derivative is positive, thus 
$h(t)<t_{\ast}$. Let $\bar{t}$ be such a point. As $r(t) \geq 0$ and 
$[t_{\ast},\tau] \subset [h(\bar{t}), \bar{t}]$, then 
$$ \int_{t_{\ast}}^{\tau} r(s)~ds \leq \int_{h(\bar{t})}^{\bar{t}}
r(s)~ds < \frac{\lambda}{L+1},
$$
consequently,
%so by (\ref{cup})
$$x(\tau)-x(t_{\ast})= \int_{t_{\ast}}^{\tau} \dot{x} (s)~ds
\leq  \int_{t_{\ast}}^{\tau} \alpha(L+1)r(s)~ds< \alpha(L+1) 
\frac{\lambda}{L+1} = \lambda \alpha.$$
Hence $x(\tau) < x(t_{\ast})+\lambda \alpha=K+\lambda \alpha$, which 
contradicts 
the assumption $x(\tau) \geq K+\lambda \alpha$. If we denote 
$t_1=t_{\ast}$ then $x(t)<K+\lambda \alpha$ for any $t \geq t_1$. 

Similar to the previous case, we can prove that there exists $\tau_1$ such 
that $x(t)\geq K-\lambda \alpha$ for $t \geq \tau_1$.

Thus, we accept $K-\lambda \alpha \leq x(t) \leq  K+\lambda \alpha$ as new 
solution bounds and consider $\tau_2$ such that $h(t) \geq \max \{ 
t_1,\tau_1 \}$ for $t
\geq \tau_2$. We repeat the induction step and obtain that there exists
$t_2 \geq \tau_2$ such that
$$ K-\lambda^2 \alpha \leq x(t) \leq K+\lambda^2 \alpha, ~~t \geq t_2,
$$
and a sequence of points $t_1 \leq t_2 \leq t_3 \leq t_n \leq \cdots$
such that
$$ K-\lambda^3 \alpha \leq x(t)\leq K+\lambda^3 \alpha, ~~t \geq t_3,
\cdots , K-\lambda^n \alpha \leq x(t) \leq K+\lambda^n \alpha, 
~~t \geq t_n, \cdots .$$
Since $0<\lambda<1$ then ${\displaystyle \lim_{t\to\infty} x(t)=K}$,
which completes the proof.
\qed

\begin{corol}
Suppose (a1)-(a2),(a4)-(a5) hold, $r>0$ and ${\displaystyle 
\sup_{t \geq 0} (t-h(t)) < \frac{1}{r(L+1)} }$ (we recall that 
$L$ is the Lipschitz constant defined in (a1)). Then all solutions of the 
equation \beq{1aadd}
\dot{x}(t) =  r \left[ \int_{h(t)}^t f( x(s)) d_s R(t,s) - x(t)
\right]
\eeq
tend to $K$.
\end{corol}

\begin{remark}
\label{Lstar}
Let us remark that the Lipschitz constant $L$ in \rf{cup} can be 
substituted by, generally, a smaller value of
\begin{equation}
\label{cup1}
L^{\ast} = \sup_{x \geq 0} \left| \frac{f(x)-K}{x-K} \right|,
\end{equation}
since in \rf{star31} only $|f(x)-K|$ is estimated for $|x-K| < \alpha$.
\end{remark}

Next, let us estimate the attracting set when (\ref{1a}) is not absolutely 
stable. We introduce $m,M$ as in \rf{mM} and
denote by $x_{\max}$ the greatest point in $[a,K]$  
where ${\displaystyle M = \max_{x \in [0,K]} f(x)}$ is attained,
by $x_{\min}$ the minimal point where ${\displaystyle m= \min_{x \in 
[K,M]} f(x)}$  is attained.

\begin{figure}[ht]   
\centering
\includegraphics[scale=0.7]{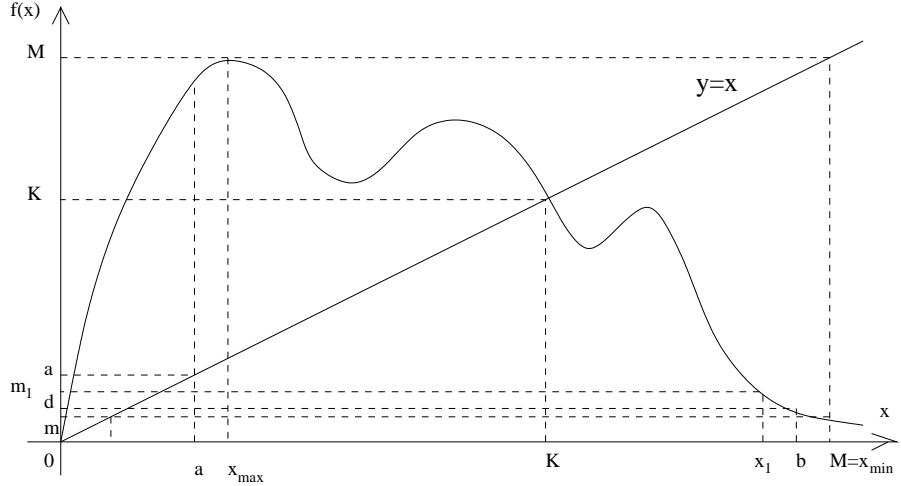}
\caption{Assuming $a$ and $b$ as in the figure exist, we can construct a delay 
equation such that $[a,b]$ is its attracting set.
}
\label{figure3}
\end{figure}

\begin{guess}
\label{theorem6}
Suppose $f(x)$ satisfies (a1) and 
\beq{2cycle}
%f(x_{\max})>x_{\min}>K, ~~~f(x_{\min})<x_{\max}<K.
m=f(x_{\min})<x_{\max}<K.
\eeq
Then for any $a\in (m,x_{\max})$, and any $b \in (K,M)$ such that 
${\displaystyle \min_{x \in [K,b]} f(x) < a}$ 
there exists a problem \rf{1a},\rf{2star} with parameters satisfying 
(a2)-(a5) such that 
\beq{attracting}
\liminf_{t \to \infty} x(t)=a, ~~\limsup_{t \to \infty} x(t)=b.
\eeq
\end{guess}
{\bf Proof.}
Let us fix $a\in (m,x_{\max})$, and $b \in (K,M)$ such that   
${\displaystyle d= \min_{x \in [K,b]} f(x) < a}$, see Fig.~\ref{figure3}.

Then there exists $x_1$ such that $K<x_1<b$ and $m_1=f(x_1)<a$.
For the initial function
%$a,b$ satisfying \rf{attracting} and set 
$$\varphi(t)=a+(t+1)(b-a),~~t \in [-1,0],$$
we have $\varphi(-1)=a$, $\varphi(0)=b$ and $\varphi$ changes continuously
from $a$ to $b$. Since $a<K<x_1<b$, then $x(s_0)=x_1$ for 
some $s_0 \in (-1,0)$. We consider the equation
\beq{particular}
\dot{x}(t)= r\left[ f(x(h(t)))-x(t) \right], ~~t \geq 0,
\eeq
where $r>0$. Obviously parameters of \rf{particular} and the initial 
function satisfy (a3)-(a5).

We choose $h(t)=s_0$, $t \in [0, \tau_1]$, where
${\displaystyle \tau_1= \frac{1}{r} \ln \left( \frac{b-m}{a-m} \right) 
>0}$. Then 
$$x(t)=(b-m)e^{-rt}+m, ~~t \in [0, \tau_1],$$ 
so $x(\tau_1)=a$ and $x(t)$ changes continuously
from $b$ to $a$ in $[0, \tau_1]$, $a<x_{\max}<b$, thus there exists $s_1 
\in (0,\tau_1)$ such that $x(s_1)=x_{\max}$. We assume 
$h(t)=s_1$, $t \in (\tau_1,\tau_2]$, where ${\displaystyle \tau_2= \tau_1
+ \frac{1}{r} \ln \left( \frac{M-a}{M-b} \right) > \tau_1 }$.
The solution is
$$x(t)= (a-M) e^{-r(t-\tau_1)} + M, ~~ t \in [\tau_1,\tau_2], $$
hence $x(\tau_2)=b$. We continue periodically with $h(t)$ piecewise 
constant such that
$$ x(h(t))=\left\{  \begin{array}{ll} x_1, & t \in (\tau_{2n}, 
\tau_{2n+1}], \\ x_{\max} & t \in (\tau_{2n+1},\tau_{2n+2} ], \end{array} 
\right.  ~~
n=0,1,2 \cdots ,
$$ 
where $\tau_0=0$, $h(t) \in (\tau_{k-1}, \tau_k)$ if $t \in (\tau_k, 
\tau_{k+1})$, $x(\tau_{2n})=b$, $x(\tau_{2n+1})=a$, $n=0,1,2 \cdots $ and
$a \leq x(t) \leq b$ for any $t$. Here the delay is bounded and piecewise
continuous, it obviously satisfies (a2).
\qed

\begin{corol}
\label{attract}
If (a1)-(a5) and \rf{2cycle} hold then the sharpest attracting interval 
for \rf{1a},\rf{2star} is $(m,M)$, where $m$ and $M$ are defined in 
\rf{mM}.
\end{corol}

\begin{remark}
\label{remark2}
Suppose the difference equation $x_{n+1}=f(x_n)$ has a stable 2-cycle 
$(a,b)$ in such a way that $f(f(x))<x$ for $a<x<K$, $x>b$ and
$f(f(x))>x$ for $x<a$, $K<x<b$. Then for any $\varepsilon>0$ there exist 
such distributed delay and initial function that 
$[a+\varepsilon,b-\varepsilon]$ is an attracting interval of 
\rf{1a},\rf{2star}. 
Moreover, if $K \in (a,b)$, function $f$ satisfies (a1), $f(a)>b$ 
and $f(b)<a$, then, similar to the proof of Theorem \ref{theorem6}
we can demonstrate that there exists problem \rf{1a},\rf{2star} such that 
$(a,b)$ is its attracting interval.
\end{remark}

\begin{remark}
\label{remark3}
The claims of Corollary \ref{attract} and Remark \ref{remark2} to some 
extent complement
Theorem 6 in \cite{RostLiz} where the sharpest invariant and attracting 
interval was 
found in the case when absolute stability fails but there exists a unique 
globally attractive 2-cycle for the relevant difference equation (for 
equations with a unique constant concentrated delay).
\end{remark}

%\begin{remark}
%Following the pattern of the proof of Theorem \ref{theorem6} we can 
%construct an example where $x(\tau_{2n}) = M-\varepsilon/2^n$,
%$x(\tau_{2n+1})=m+\varepsilon/2^n$, so the attracting set will be exactly
%$m,M$ but the length of segments $\tau_k, \tau_{k+1}$ will be growing.
%Here the delay will not be bounded anymore but (a2) is still satisfied.
%\end{remark}

%\begin{remark}
%\label{remark5}
%In the assumptions of Theorem \ref{theorem6} constant $L$ can be
%substituted by, generally, a smaller constant
%\begin{equation}
%\label{properL}
%L^{\ast} = \sup_{x>0}\left| \frac{f(x)-f(K)}{x-K} \right|.
%\end{equation}
%\end{remark}

\section{Applications and Discussion}

The global attractivity results can be applied to any 
unimodal function $f$ and any type of finite delay.
For example, we can deduce the following result in the stream of 
\cite{RostWu}, see also \cite{Singer} and Proposition 2.1 in \cite{Liz}.

\begin{guess}
\label{maintheoremapplic}
Suppose (a1)-(a5) hold, $f$ is three times continuously differentiable and has the 
only critical point $x_0>0$ (maximum),
$$ (Sf)(x)= \frac{f^{\prime\prime\prime}(x)}{f^{\prime} (x)} -
\frac{3}{2} \left( \frac{f^{\prime\prime}(x)}{f^{\prime} (x)} \right)^2 
<0 \mbox{  ~  for ~~}x \neq x_0 $$
and $|f^{\prime}(K)| \leq 1$.
%$f^{\prime}(x)<0$ if $x \geq x_0$, $f^{\prime\prime}(x)<0$ if $0 \leq x 
%\leq x_0$, $\lim\limits_{x \to \infty} f(x)=0$ and 
%\beq{Lcondition}
%\alpha= f(f(x_0))>x_0.
%\eeq
Then any solution of (\ref{1a}),\rf{2star} tends to $K$ as $t \to \infty$.
\end{guess}

As an application, consider the Nicholson's blowflies equation with a 
distributed delay
\beq{71new}
\dot{x}(t)= -\delta x(t)+  p\int_{h(t)}^t x(s) e^{-a x(s)}~d_s R(t,s),
\eeq
where (a2)-(a5) hold.

If $p<\delta$ then by Theorem \ref{maintheorem2} all positive solutions 
of \rf{71new} go to extinction, i.e., $x(t) \to 0$ as $t \to \infty$.
If ${\displaystyle 1 < \frac{p}{\delta} < e^2}$ then by Theorem 
\ref{maintheorem1} all positive 
solutions of \rf{71new} tend to the positive equilibrium ${\displaystyle 
K=\frac{1}{a} \ln \left( \frac{p}{\delta} \right)}$.
However, since sustainable oscillations were observed in experiments
\cite{nichol}, we speculate that in the real problem the ratio $p/\delta$
was above $e^2$ and delays significant enough.

Particular cases of \rf{71new} are the equation with a variable delay
\beq{71newa}
\dot{x}(t)= -\delta x(t)+  p x(h(t)) e^{-a x(h(t))},
\eeq
including the original equation \rf{71} with a constant delay,
the integrodifferential equation
\beq{71newb}
\dot{x}(t)= -\delta x(t)+  p \int_{h(t)}^t k(t-s) x(s) e^{-a x(s)} \, ds,
\eeq
where
\beq{kernel}
\int_{h(t)}^t k(t-s) ~ds =1 \mbox{~ for any ~} t \geq 0,
\eeq
and the mixed equation
\beq{71newc}
\dot{x}(t)= -\delta x(t)+  \alpha p  x(g(t)) e^{-a x(g(t))}
+ (1-\alpha) p \int_{h(t)}^t k(t-s) x(s) e^{-a x(s)} \,ds,
~~0 \leq \alpha \leq 1.
\eeq

By Theorems \ref{maintheorem1}, \ref{maintheorem2} and \ref{theorem5} we 
conclude the following result.
\begin{guess}
\label{nicholson}
If
$p \leq \delta$ then all positive solutions of \rf{71new} (and thus any 
of \rf{71newa},\rf{71newb},\rf{71newc} ) go to 
extinction, if $\delta < p < \delta e^2$, all positive solutions
tend to the positive equilibrium.
If $p \geq e^2$ and 
\begin{equation}
\label{nick}
\limsup_{t \to \infty} (t-h(t))< \frac{1}{p+\delta},
\end{equation}
then all positive solutions of \rf{71new} tend to the positive 
equilibrium.
\end{guess}
{\bf Proof.} The results for $p < \delta e^2$ are immediate corollaries
of Theorems \ref{maintheorem1} and \ref{maintheorem2}. Since the 
absolute value of the derivative of ${\displaystyle f(x)= 
\frac{p}{\delta}x e^{-ax}}$ does not 
exceed $p/\delta$ (which is attained at $x=0$), then \rf{cup} implies 
\rf{nick}.
\qed

\begin{remark}
\label{Lstar1}
We can apply Remark \ref{Lstar} to improve condition
(\ref{nick}). In fact, $L^{\ast}$ defined in \rf{cup1} 
does not exceed $p/(\delta e^2)$ (the absolute value of the 
minimum of the derivative which is attained at $x=2/a$)
if $p \geq \delta e^2$ and the stability condition becomes
\begin{equation}
\label{nick1}
\limsup_{t \to \infty} (t-h(t)) < \frac{1}{p/e^2+\delta}.
\end{equation}
We also notice that for unimodal functions such that their equilibrium 
(the fixed point $x=f(x)$) exceeds maximum point $x_{\max}$ and such that the 
minimum of the derivative is less than -1 (otherwise, under some 
additional conditions which are satisfied for the Nicholson's blowflies 
and the Mackey-Glass equations, the relevant difference equations are 
stable \cite{Liz}), then $L^{\ast}$ defined in \rf{cup1} does not exceed
the absolute value of the minimum of the derivative. 

Really, function ${\displaystyle g(x) = \left| \frac{f(x)-K}{x-K} 
\right|}$ satisfies $g(0)=1$, ${\displaystyle 0< g(x) \leq - 
\!\! \min_{x \in  [x_{\max}, \infty )} \!\!\! f^{\prime} (x) }$ in 
$[x_{\max}, \infty )$. The value in the right 
hand side of the latter inequality exceeds 1 and $g$ is decreasing in $[0, 
x^{\ast}]$, where $x^{\ast}<x_{\max}$ is such point that $f(x^{\ast})=K$, while in 
$[x^{\ast},x_{\max})$ function $g(x)$ is less than $g(x_{\max})$, which
does not exceed
the absolute value of the minimum of the derivative.
\end{remark}

Let us compare Theorem \ref{nicholson} and Remark \ref{Lstar1} to 
some known global stability results 
for the Nicholson's blowflies equation. For example, \cite{GL}
contains the following global stability condition
\beq{GLcond}
\left( e^{\delta \tau}-1 \right) \left( \frac{p}{\delta} -1 \right)<1
\eeq
for equation \rf{71} with a constant delay.

\begin{figure}[ht]
\centering
\includegraphics[scale=0.5]{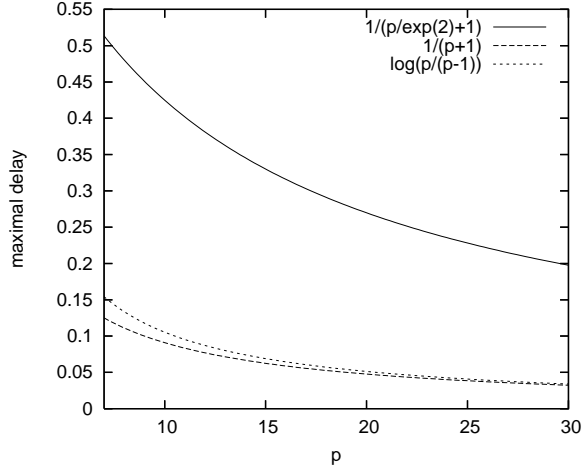}
\caption{The values of delay for which
(\protect{\ref{71}}) is asymptotically stable are given for $\delta 
=1$, $p \geq \delta e^2 >7$ and conditions (\protect{\ref{nick1}}), 
(\protect{\ref{nick}}) and (\protect{\ref{GLcond}}), respectively.
}
\label{fignick}
\end{figure}

Condition \rf{GLcond} is slightly sharper for \rf{71} than \rf{nick} but 
\rf{nick} is applicable to equations with variable delays and 
integrodifferential equations as well. Condition \rf{nick1} is sharper 
than \rf{GLcond}, see Fig. \ref{fignick}. It is also applicable to 
equations with a variable and/or distributed delay.

% plot [7:30] 1/(x/exp(2)+1),1/(x+1),log(x/(x-1))
% set xlabel "p"
% set ylabel "maximal delay"
% set term postscript portrait
% set size 0.9,0.5
% set output "figure_nickol.ps"
% replot

For completeness of references, we remark that oscillation of \rf{71new}
was studied in \cite{MCM2008}.

All the main results are also relevant for the Mackey-Glass equation with
a distributed delay
\beq{72new}
\dot{x}(t) = \int_{h(t)}^t \frac{ax(s)}{1+x^{\gamma} (s)}~d_s 
R(t,s)-bx(t),~a>0,~b>0, ~\gamma>0,
\eeq
which also involves the equation
\beq{72newc}
\dot{x}(t)= -b x(t)+  \frac{\alpha a  x(g(t))}{1+x^{\gamma} (g(t))}  
+ (1-\alpha) a \int_{h(t)}^t k(t-s) \frac{x(s)}{1+x^{\gamma} (s)}~ds,
~~0 \leq \alpha \leq 1,
\eeq
where \rf{kernel} is satisfied,
as well as its special cases of the equation with a variable delay 
($\alpha=1$) and the integrodifferential equation ($\alpha=0$).
If $a<b$ then all solutions of \rf{72new} go to extinction. If $a>b$ and
\beq{condMackey}
\gamma < \frac{2a}{a-b}
\eeq
then all solutions tend to the positive equilibrium ${\displaystyle
K=\left( \frac{a}{b} -1 \right)^{1/\gamma} }$.

Again,  Theorems \ref{maintheorem1}, \ref{maintheorem2} and 
\ref{theorem5} imply the following result.
\begin{guess}
\label{mackeyglass}
If
$a \leq b$ then all positive solutions of \rf{72new} 
(and thus of \rf{72newc} ) go to 
extinction. If either $0< \gamma \leq 2$ or $\gamma>2$ and ${\displaystyle 
b < a < \frac{\gamma b}{\gamma-2} }$, then all positive solutions
tend to the positive equilibrium $K$.
%${\displaystyle x^{\ast} = \left( \frac{a-b}{b} \right)^{1/\gamma} }$.
If $\gamma >2$, ${\displaystyle a \geq \frac{\gamma b}{\gamma-2} }$
and 
\begin{equation}
\label{mack}
\limsup_{t \to \infty} (t-h(t))< \frac{4b^2\gamma}{a^2(\gamma-1)^2 +4 a b 
\gamma},
\end{equation}
then all positive solutions of \rf{72new} tend to the positive 
equilibrium.
\end{guess}
{\bf Proof. } 
We remark that for $\gamma \leq 1$ the function ${\displaystyle f(x)
= \frac{ax}{b(1+x^{\gamma})}}$ is monotone increasing, and we have 
stability for any delay. For $1 \leq \gamma \leq 2$ the relevant 
difference equation $x_{n+1}=f(x_n)$ is locally asymptotically 
stable, we apply Theorem \ref{maintheoremapplic}.
It is known that for the relevant difference 
equation local asymptotic stability implies global asymptotic stability
(see \cite{Liz}, Example 3.2), thus the difference equation is globally 
asymptotically stable whenever
$$
|f^{\prime} (K)|= \left| 1-\gamma+ \frac{b}{a} \gamma \right|<1,
~~\mbox{ where ~~} \gamma > 2,$$
or $a< \gamma b/(\gamma-2)$. According to Remark \ref{Lstar1}, we can 
take the absolute value of the minimum of the derivative which is attained at 
${\displaystyle \left( \frac{\gamma+1}{\gamma-1}\right)^{1/\gamma} }$
as $L^{\ast}$:
$$
L^{\ast}=\frac{a(\gamma-1)^2}{4b\gamma}.
$$
Application of Theorem \ref{theorem5} and Remark
\ref{Lstar1} completes the proof.
\qed

Finally, let us formulate some open problems.

\begin{enumerate}
\item
Generally, absolute stability results are incorrect for systems with 
infinite memory. Deduce stability results when (a2) is not satisfied
but measure $d_s R(t,s)$ decays exponentially with memory
$$
\int_{-\infty}^{t-\tau} d_s R(t,s) \leq S e^{-\alpha \tau}.
$$
\item
For arbitrary $f(x)$ satisfying (a1) and the equation with a constant 
delay, characterize attractive sets.
\item
Obtain stability results for equations with a distributed delay in the 
multistability case: there are several points satisfying $f(x)=x$.
If $K_1$ is the first positive fixed point and $f(x)<x$ for $0<x<K_1$,
describe initial conditions leading to extinction (Alley effect). 
\item
Would the results of the present paper remain valid for the equation
$$
\dot{x}(t) =  r(t) \left[ f\left( \int_{h(t)}^t x(s) d_s R(t,s) \right) - x(t)
\right], ~t \geq  0,
$$
if (a1)-(a5) hold?
\end{enumerate}


\begin{thebibliography}{99}
\bibitem{KWW1999}
T. Krisztin, H.\,O. Walther and J. Wu, Shape, smoothness and 
invariant stratification of an attracting set for delayed monotone 
positive feedback. Fields Institute Monographs, 11. American Mathematical 
Society, Providence, RI, 1999.
\bibitem{KW2001}
T. Krisztin and H.\,O. Walther, Unique periodic orbits for 
delayed positive feedback and the global attractor, 
{\em  J. Dynam. Differential Equations} {\bf 13}  (2001),  no. 1, 1--57.  
\bibitem{MPS}
J. Mallet-Paret and G.\,R. Sell, The Poincar\'{e}-Bendixson 
theorem for monotone cyclic feedback systems with delay, {\em J. Differential 
Equations} {\bf 125} (1996), 441–-489.
\bibitem{MG1}
M.\,C. Mackey and L. Glass, Oscillation and chaos in physiological control
systems, {\em Science} {\bf 197} (1977), 287--289.
\bibitem{LW}
B. Lani-Wayda, Erratic solutions of simple delay equations, {\em 
Trans. Am. Math. Soc.} {\bf 351} (1999), 901???945. 
%(doi:10.1090/S0002-9947-99-02351-X). 
\bibitem{RostWu}
G. R\"{o}st and J. Wu, Domain-decomposition method for the 
global dynamics of delay differential equations with unimodal feedback,  
{\em Proc. R. Soc. Lond. Ser. A Math. Phys. Eng. Sci.} {\bf 463}  (2007),  
no. 2086, 2655--2669. 
\bibitem{RostLiz}
E. Liz and G. R\"{o}st, On the global attractor of delay differential 
equations with unimodal feedback, submitted.
\bibitem{Singer}
D. Singer, Stable orbits and bifurcation of maps of the interval,  
{\em SIAM J. Appl. Math.} {\bf 35}  (1978), no. 2, 260--267. 
\bibitem{IvShark}
A.\,F. Ivanov and A.\,N. Sharkovsky, Oscillations in 
singularly perturbed delay equations, {\em Dynam. Rep. (New Series)} {\bf 
1} (1992), 164–-224.
\bibitem{nichol}
A.\,J. Nicholson, An outline of the dynamics of animal populations,
{\em Austral. J. Zool.} {\bf 2} (1954), 9--65.
\bibitem{GBN}
W.\,S.\,C. Gurney, S.\,P. Blythe and R.\,M. Nisbet,
Nicholson's blowflies revisited, {\em Nature} {\bf 287} (1980), 17--21.
\bibitem{CAMQ2006}
E. Braverman and D. Kinzebulatov, Nicholson's blowflies equation
with a distributed delay,
{\em Can. Appl. Math. Quart.} {\bf 14} (2006), no. 2, 107--128.
\bibitem{MG2}
J. Losson, M.\,C. Mackey, and A. Longtin, Solution multistability in first
order nonlinear differential delay equations, {\em Chaos} {\bf 3} (1993), no. 2, 
167--176.
\bibitem{Liz1}
E. Liz, E. Trofimchuk and S. Trofimchuk, Mackey-Glass type
delay differential equations near the boundary of absolute stability,
{\em  J. Math. Anal. Appl.} {\bf 275}  (2002),  no. 2, 747--760.
\bibitem{Liz2}
E. Liz, C. Mart\'{i}nez and S. Trofimchuk, Attractivity
properties of infinite delay Mackey-Glass type equations, {\em
Differential Integral Equations} {\bf 15}  (2002),  no. 7, 875--896.
\bibitem{KubSak}
I. Kubiaczyk and S.\,H. Saker, Oscillation and stability in nonlinear
delay differential equations of population dynamics, {\em Math. Comput.
Modelling} {\bf 35}  (2002),  no. 3-4, 295--301.
\bibitem{Saker}
S.\,H. Saker, Oscillation and global attractivity in hematopoiesis
model with delay time, {\em Appl. Math. Comput.} {\bf 136}  (2003),  no.
2-3, 241--250.
\bibitem{BB2006}
L. Berezansky and E. Braverman,
Mackey-Glass equation with variable coefficients, {\em Comput. Math.
Appl.} {\bf 51}  (2006),  no. 1, 1-16.
\bibitem{BBDCDIS}
L. Berezansky and E. Braverman, On existence and attractivity of
periodic solutions for the hematopoiesis equation,
{\em Dyn. Contin. Discrete
Impuls. Syst. Ser. A Math. Anal.} {\bf 13B}  (2006),  suppl., 103--116.
\bibitem{Shaikhet}
N. Bradul and L. Shaikhet, Stability of the positive point of 
equilibrium of Nicholson's blowflies equation with stochastic 
perturbations: numerical analysis, {\em  Discrete Dyn. Nat. Soc.} {\bf  
2007}, Art. ID 92959, 25 pp. 
\bibitem{Cordun}
C. Corduneanu, Functional Equations with Causal Operators. {\em Stability
and Control: Theory, Methods and Applications}, {\bf 16}. Taylor \&
Francis, London, 2002.
\bibitem{GL} 
I. Gy\"{o}ri and G. Ladas,
Oscillation Theory of Delay Differential Equations,
Clarendon Press, Oxford, 1991.
\bibitem{Li}
J. Wei and M.\,Y. Li, Hopf bifurcation analysis in a
delayed Nicholson blowflies equation, {\em Nonlinear Anal.} {\bf 60}
(2005),  no. 7, 1351--1367.
\bibitem{Coppel}
W.\,A. Coppel, The solution of equations by iteration, {\em Proc.
Cambridge Philos. Soc.} {\bf 51} (1955), 41--43.
\bibitem{Liz}
E. Liz,  Local stability implies global stability in
some one-dimensional discrete single-species models,
{\em Discrete Contin. Dyn. Syst. Ser. B}~ {\bf 7} (2007), no. 1, 191--199.
\bibitem{MCM2008}
L. Berezansky and E. Braverman, Linearized oscillation theory for a 
nonlinear equation with a distributed delay,
{\em Math. Comput. Modelling} {\bf 48} (2008), no. 1-2, 
287-304.


\end{thebibliography}
\end{document}